\documentclass[12pt]{article}
\usepackage{epsfig}
\usepackage{amsmath}

\newcommand{\no}{\nonumber}

\newcommand{\be}{\begin{equation}}

\newcommand{\ee}{\end{equation}}

\newcommand{\ba}{\begin{eqnarray}}

\newcommand{\ea}{\end{eqnarray}}

\begin{document}
\title{\hspace{.2 in} \\ 
{\Large \bf An Orthogonal Discrete Auditory Transform}}

\author{Jack Xin \thanks{Corresponding author, Department of Mathematics and ICES,
University of Texas at Austin,
Austin, TX 78712, USA; email:jxin@math.utexas.edu.} \hspace{.1 in} and 
\hspace{.1 in} Yingyong Qi \thanks{
Qualcomm Inc,
5775 Morehouse Drive,
San Diego, CA 92121, USA.}
}

\date{}

\maketitle

\vspace{0.25 in}
\thispagestyle{empty}

\begin{abstract}
An orthogonal discrete auditory transform (ODAT) from sound signal
to spectrum is constructed by combining the auditory spreading matrix of Schroeder et al 
and the time one map of a discrete nonlocal Schr\"odinger equation.
Thanks to the dispersive smoothing property of the Schr\"odinger evolution, 
ODAT spectrum is smoother than that of  
the discrete Fourier transform (DFT) consistent with human audition. 
ODAT and DFT are compared in signal denoising tests with spectral thresholding method. The signals 
are noisy speech segments. ODAT outperforms DFT in signal to noise 
ratio (SNR) when the noise level is relatively high. 
\end{abstract}

\vspace{.4 in}

\hspace{.2 in} {\bf Keywords: Orthogonal Discrete
Auditory Transform,}
\vspace{ .1 in} 

\hspace{1.2 in} {\bf Schr\"odinger Equation.}

\vspace{.4 in}

\hspace{.2 in} {\bf AMS Subject Classification: 94A12, 94A14, 65T99.}
            
\newpage
\setcounter{page}{1}

\section{Introduction}
\setcounter{equation}{0}
Acoustic signal processing can benefit significantly from utilizing properties of human 
audition, e.g. perceptual coding in MP3 technology of music compression \cite{Poh,S}. 
In \cite{xq05}, an invertible discrete auditory transform (DAT) is formulated by the present 
authors to map sound signal to auditory spectrum. DAT is more adapted to the spectral 
features of the ear than Fourier transform. It incorporates the auditory 
spreading functions of Schroeder, Atal and Hall \cite{S} to achieve 
smoother spectrum than that of the discrete Fourier transform (DFT), and better 
performance in denoising under spectral thresholding. However, such a transform has redundancy 
in the sense that the image of a discrete vector lies in a higher dimensional space, similar to 
tight frames in wavelets \cite{ID,Strang}.   
\medskip

In this paper, such redundancy is removed by constructing an orthogonal (unitary) matrix with 
spreading property over frequency bands comparable to the critical bands in hearing. 
Critical bands (Table 10.1, p 309, \cite{Poh}) characterize the bandwidth of the human auditory filter.
The auditory orthogonal matrix is obtained from the time one map of a spatially discrete 
nonlocal Schr\"odinger equation. The Schr\"odinger equation conserves the $L^2$ norm or 
Euclidean length, implying the orthogonality of the time one map. On the other hand, 
the dispersive smoothing nature of the Schr\"odinger evolution leads to 
the spreading property of the time one map.
The auditory functions of Schroeder, Atal and Hall \cite{S} appear as a nonlocal potential 
in the Schr\"odinger equation. As a result, a class of orthogonal discrete 
auditory transforms (ODAT) are generated. In searching for ODATs, an alternative method 
based on the dilation equation of wavelets is also found, however, such an approach turns out to be 
too rigid to accomodate auditory properties, e.g. spectral spreading across 
critical bands. 
\medskip

The paper is organized as follows. In section 2, the ODAT is derived from the general 
DAT \cite{xq05}, and the ODAT construction is presented based on  
the discrete Schr\"odinger equation. A specific ODAT is given by inserting the auditory 
spreading functions in \cite{S}.
In section 3, auditory spectra of a two tone signal (with frequencies across a critical band) and 
of a vowel segment are compared with their DFT counterparts to illustrate the auditory 
spectral spreading. Denoising with spectral thresholding is performed on voiced and unvoiced speech 
segments. ODAT is found to increase signal to noise ratio beyond DFT when 
the noise content is relatively high.
Concluding remarks are made in section 4. 

\section{ODAT and Schr\"odinger}
\setcounter{equation}{0}
Let $s=(s_0,\cdots,s_{N-1})$ be a discrete real signal, 
the discrete Fourier transform (DFT) is \cite{Br}:
\be
 \hat{s}_{k} = \sum_{n=0}^{N-1}\, s_{n}\, e^{-i (2\pi n k/N )}. \label{dt1}
\ee
The general discrete auditory transform (DAT) is \cite{xq05}:
\be
S_{j,m} \equiv \sum_{l=0}^{N-1}\, s_l \, K_{j-l,m}, \label{dt5}
\ee
where the double indexed kernel function is:
\be
K_{l,m} = \sum_{n=0}^{N-1}\, X_{m,n} \, 
e^{i (2\pi l n/N)}; \label{dt6}
\ee
where the matrix $X_{m,n}$ has square sum equal to one in $m$:
\be
\sum_{m=0}^{M-1}\, |X_{m,n}|^{2}  =1, \;\; \forall n. 
\label{nm}
\ee
Here $M$ is on the order of $N$.
\medskip

DFT is recovered from DAT by setting $j=0$, $M=N$, and $X_{m,n}$ the 
$N\times N$ identity matrix.  
In case that $X_{m,n}$ is a nontrivial orthogonal matrix, let us still set $j=0$ 
in (\ref{dt5}) to find:
\ba
S_{0,m} \equiv S_m & = & \sum_{l=0}^{N-1}\, s_l\, \sum_{n=0}^{N-1}\, X_{m,n} \, 
e^{- i (2\pi l n/N)} \no \\
& = & \sum_{n=0}^{N-1}\, X_{m,n}\, \hat{s}_{n}. \label{dt7}
\ea
The mapping from $s_{l}$ to $S_m$ is orthogonal. The problem reduces to finding  
an orthogonal matrix $(X_{m,n})$ with auditory features. 
\medskip

Such a matrix acts on complex numbers $\hat{s}_{n}$ (except the modes $n=0$ and 
$n=N/2$, so called DC and Nyquist modes). Let us consider the time one map of 
the following spatially discrete Schr\"odinger equation:
\be
i\, u_{n,t} = \sigma_1 \, (u_{n+1} - 2\, u_{n} + u_{n-1}) + 
\sigma_2\, \sum_{m=1}^{N_h} \, V_{m,n}\, u_{m}, \label{dt8}
\ee
where $\sigma_1$ and $\sigma_2$ are positive real numbers, $N_h = N/2 -1$, $(V_{m,n})$ is 
a symmetric $N_h \times N_h$ matrix to carry certain auditory information of the ear. 
For simplicity, Dirichlet boundary condition is imposed for the evolution of 
equation (\ref{dt8}). The discrete equations (\ref{dt8}) 
can be cast in the matrix form:
\be 
i\, U_t = (\sigma_1 \, A + \sigma_2 \, B)\, U, \label{dt9}
\ee
where $U=(u_1,u_2,\cdots,u_{N_h})^{T}$, $A$ the tridiagonal matrix ($-2$ on the diagonal, 
$1$ on the two off-diagonals), $B$ the real symmetric matrix with entry $V_{m,n}$ at $(m,n)$.
The time one map of (\ref{dt9}), denoted by $T_w$, is simply $\exp\{ i (\sigma_1 \, A + \sigma_2\, B)\}$ 
which is clearly orthogonal, $T_w \, T_{w}^{'} = Id_{N_h}$, where the prime denotes the  
conjugate transpose. 
\medskip

The matrix $B$ is built from 
auditory spreading functions \cite{S} denoted by $S(b(f_m),b(f_n))$, where $f_m$ is the frequency 
to spread from, $f_n$ is the frequency to spread to, and $b$ is the standard mapping from 
Hertz (Hz) to Bark scale \cite{Hart}. The functional form of $S(\cdot,\cdot)$ is given in 
\cite{S}. Define $V_{m,n}= 1/2\cdot (S(b(f_m),b(f_n)) + S(b(f_n),b(f_m)))$, so $B$ is the symmetric 
part of the matrix $(S(b(f_m),b(f_n))$. Numerical results based on this choice of $B$ will be reported 
in the next section.
\medskip

The matrix $X=(X_{m,n})$ takes the block diagonal form:
\be X = {\rm diag} \{ 1, T_w, 1, \widehat{T_w}^{*} \}, \label{dt10}
\ee
where the tilde denotes the reverse permutation of columns of $T_w$ so that 
the spreading occurs symmetrically
on the DFT components ($\hat{s}_{l}$, $N_h +2\leq l \leq N-1$) to preserve the conjugate 
symmetry of the spectrum.  The matrix $X$ is clearly orthogonal and leaves invariant the 
DC and Nyquist modes. The ODAT matrix is the product of $X$ and DFT matrix.  
\medskip

The continuum version of (\ref{dt8}) is:
\be
i\, u_t = \Delta_x \, u + V(x)*u,\;\; \;x\in R^n, \;\;\; n \geq 1, \label{dt11}
\ee
where $*$ is convolution, $V(x)$ is real and even. The $L^2$ norm of $u$ is conserved 
in time. Schr\"odinger equations analogous to (\ref{dt11}) have been much
studied regarding smoothing (scattering) properties and derivation from particle dynamics,
\cite{EYau04,Gin89,HO89,Jensen86,Kap96} among others. When the convolution term is cubically 
nonlinear in $u$, the equation is known as Schr\"odinger-Hartree \cite{EYau04,HO89,Gin89}. 
In \cite{HO89,Jensen86}, the smoothing and spreading property is measured in the weighted norm 
$\|\psi \|_{m,s} = \|(1+|x|^2)^{s/2}(1 -\Delta)^{m/2}\, \psi \|_{2}$, $\Delta$ the spatial Laplacian. 
Solutions at time $t \not = 0$ satisfy the bound:

\be \|u(t)\|_{1,-1} \leq C(\|u(0)\|_{0,1})\, (|t|+ |t|^{-1}). \label{dt12}
\ee
We shall see in the next section that the time one Schr\"odinger map $T_w$ 
inherits the smoothing and spreading property of the continuum case. 

\section{Numerical Tests}
\setcounter{equation}{0}
The computation is carried out in Matlab, with ODAT parameters  
$(\sigma_1,\sigma_2)= (0.6,0.04)$. Discrete signal (frame) length $N=256$. 
First consider a two tone signal consisting of sinusoids of 
frequencies 3 kHz (kilo-Hertz) and 4.3 kHz with identical amplitudes. The two frequency values  
span a critical band. Figure 1 compares the ODAT (dashed) and DFT (solid) log-magnitude spectra. 
The ODAT spectral peak regions are lower 
and wider than DFT's. Also there is more spreading in ODAT spectrum towards higher 
frequency, consistent with upward masking property of human ear \cite{ZF}.
This can be explained by the weighted norm estimate (\ref{dt12}), where large $x$ corresponds 
to large frequency. Figure 2 shows ODAT and DFT spectra of a vowel segment containing multiple 
harmonics, spectral smoothing is observed again. 
\medskip

ODAT and DFT were used to denoise speech signals via the  
thresholding method in the transformed domain \cite{xq05}.
The aim is to improve the signal-to-noise ratio (SNR) of noisy speech. The premise of the method 
is that low level components in the transformed domain are more likely to be noise 
than signal plus noise. So thresholding could improve the overall SNR of 
the signal. The simple thresholding method serves to illustrate 
the difference between ODAT and DFT in signal processing. 
A vowel and a consonant speech segments were selected,
each segment has 512 data points. Noisy speech was created by adding 
Gaussian noise to the selected segments. The level of noise was set to produce the 
SNR ranging from -12 decible (dB) to +12 dB with a 3 dB step size. ODAT and DFT were applied 
to the noisy speech signals. The magnitude of transformed components were then 
compared to a threshold. All components with magnitude smaller than the threshold 
were ignored for the reconstruction of the signal. The threshold was computed 
as the average of the DFT magnitude spectrum. Signal was reconstructed by 
the inverse ODAT and DFT, respectively. The SNRs of the reconstructed vowel signal 
is plotted vs. input SNRs in Figure 3. The SNRs of the reconstructed consonant signal is 
plotted vs. input SNRs in Figure 4.
We see that ODAT (solid) improves over DFT (dashdot) in terms of SNR 
when the noise level is relatively high, particularly in case of consonants which 
resemble noise more than the vowels. 
\medskip

The noise-reduction advantage can be attributed to the spectral 
spreading property of ODAT. The redundant DAT \cite{xq05} is quite similar 
in this respect. Redundancy however renders more modes in the transformed 
domain, and was observed to provide more SNR gain in denoising tests. It is 
interesting to find out how to enhance the amount of smoothing for ODAT in future work. 

\medskip

It is rewarding to investigate how well a nonlinear nonlocal 
Schr\"odinger equation can model the ear's nonlinear responses. Ear's nonlinearities are 
nonlinear and nonlocal in nature and the physiological models are dispersive nonlinear 
nonlocal \cite{JGeis83,deng93,deboernut,xq04,xqd}.

\section{Concluding Remarks}
Orthogonal discrete auditory transforms (ODAT) are introduced based on nonlocal spatially 
discrete Schr\"odinger equations. Dispersive smoothing, mass conservation, and robustness of 
the Schr\"odinger equation allows one to inject auditory knowledge in the transform while preserving 
orthogonality. Numerical tests on two tone and speech segments demonstrate the 
spectral spreading property of ODAT and advantage in denoising. Future work will explore 
efficient ways to enhance spectral spreading for ODATs and more complex 
signal processing applications.

\section{Acknowledgements}
This work was supported in part by NSF grant ITR-0219004, and 
a Fellowship from the John Simon Guggenheim Memorial Foundation. 
We thank Prof. G. Papanicolaou and Prof. H-T Yau for  
helpful conversations. Results in the paper were presented at the IPAM workshop 
on ``Mathematics of the Ear and Sound Signal Processing'', UCLA, Jan 31 -- Feb 2, 2005. 
We thank many workshop speakers and participants for their comments and interest. 

\vspace{.2 in}

\bibliographystyle{plain}








\newpage
\begin{figure}[p]
\centerline{\includegraphics[width=350pt,height=340pt]{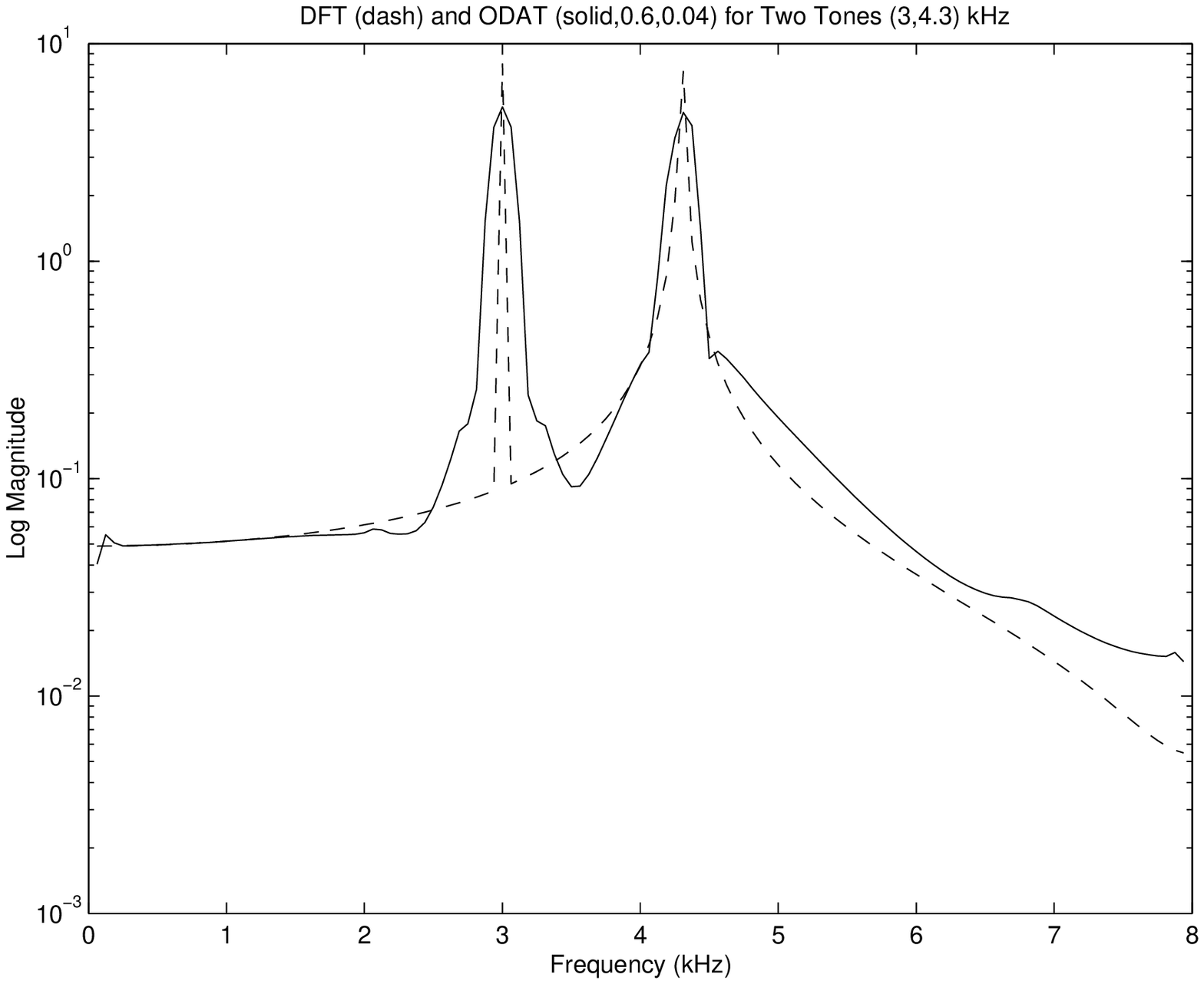}}
\vspace{.25 in}
\label{Fig1}
\caption{Comparison of ODAT spectrum (solid) and DFT spectrum (dash) of a two tone 
signal of frequencies (3,4.3) kHz and identical amplitudes. ODAT's spectral spreading appears 
near the peak areas and towards the higher frequencies. 
ODAT parameters $(\sigma_1,\sigma_2)=(0.6,0.04)$.}
\end{figure}

\newpage
\begin{figure}
\centerline{\includegraphics[width=350pt,height=350pt]{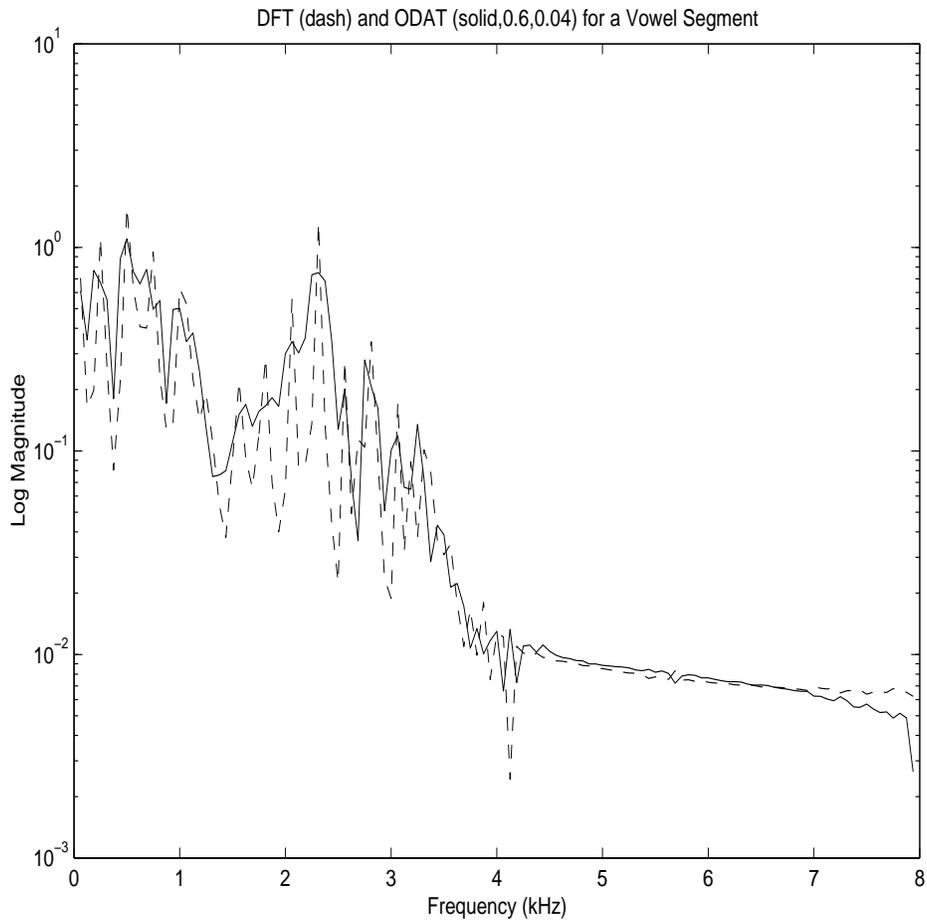}}
\vspace{.25 in}
\caption{Comparison of ODAT spectrum (solid) and DFT spectrum (dash) of a vowel segment. 
ODAT parameters $(\sigma_1,\sigma_2)=(0.6,0.04)$.}
\label{Fig2}
\end{figure}

\newpage
\begin{figure}
\centerline{\includegraphics[width=350pt,height=350pt]{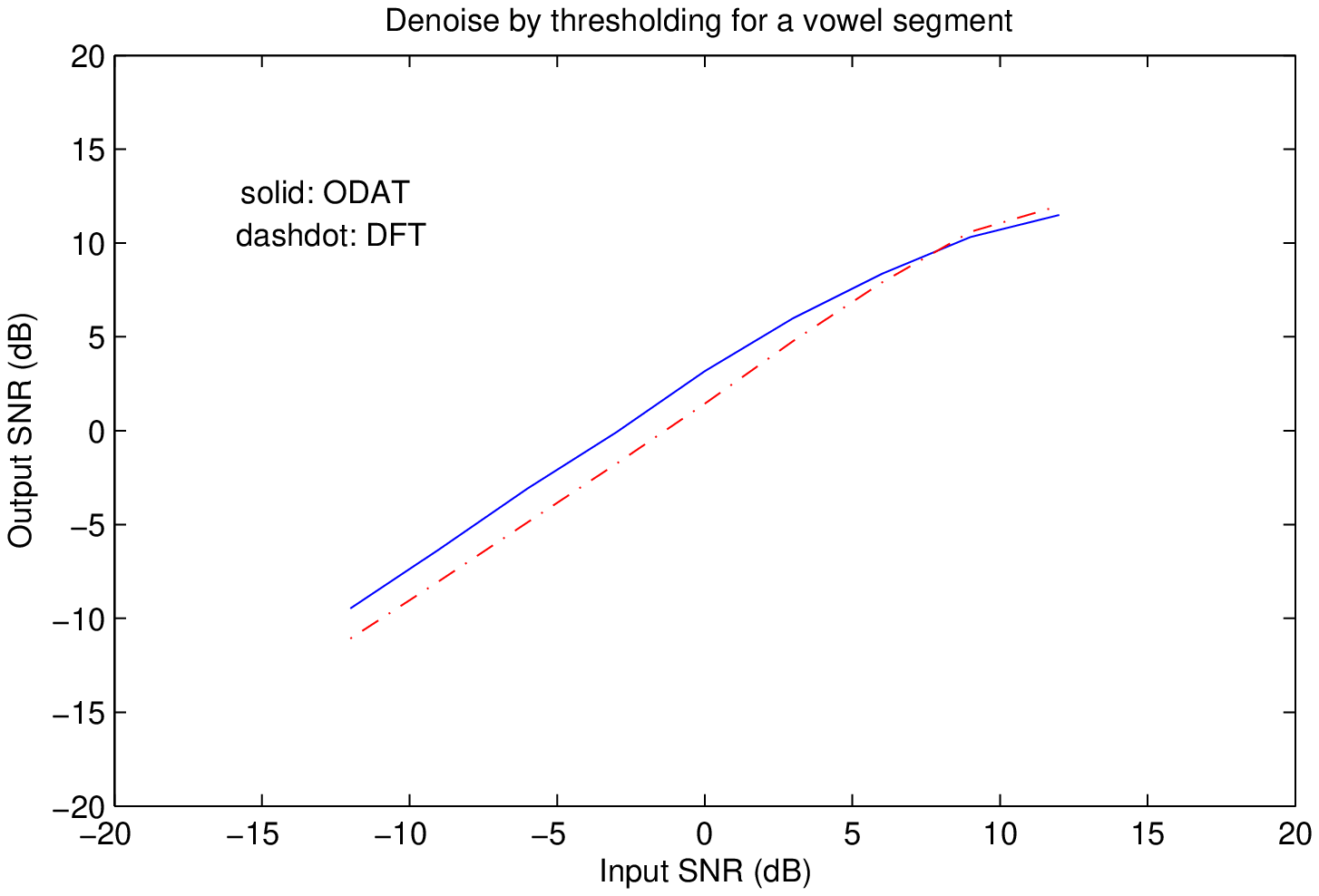}}
\vspace{.25 in}
\caption{Comparison of ODAT (solid) and DFT (dashdot) denoising by spectral thresholding for 
a vowel segment. Spectral spreading property of ODAT helps to increase signal 
content when noise level is relatively high, e.g. input SNR below 7 decible (dB).}
\label{Fig3}
\end{figure}

\newpage
\begin{figure}
\centerline{\includegraphics[width=350pt,height=350pt]{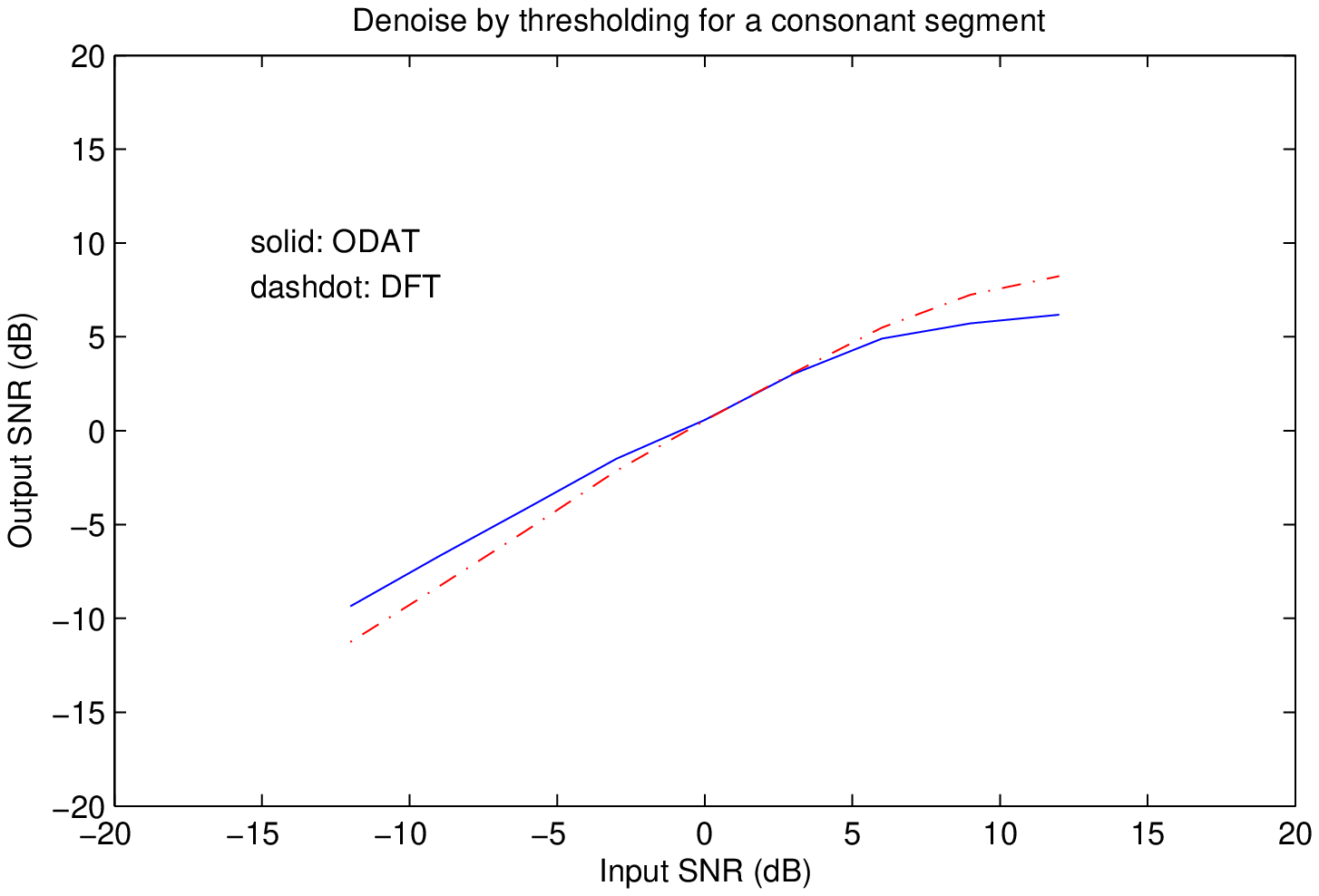}}
\vspace{.25 in}
\caption{Comparison of ODAT (solid) and DFT (dashdot) denoising by spectral thresholding for 
a vowel segment. Spectral spreading property of ODAT helps to increase signal 
content when noise level is relatively high, e.g. input SNR below zero decible (dB).}
\label{Fig4}
\end{figure}

\end{document}